\documentclass[aap,seceqn,MSNbibl,citesort,dvips]{arximspdf}

\doi{10.1214/09-AAP675}
\volume{20}
\issue{5}
\pubyear{2010}
\firstpage{1891}
\lastpage{1906}

\makeatletter

\newproclaim{asm}{Assumption}[section]
\newproclaim{dfn}[asm]{Definition}
\newtheorem{theorem}[asm]{Theorem}
\newtheorem{lem}[asm]{Lemma}
\newproclaim{rem}[asm]{Remark}

\makeatother

\begin{document}
\begin{frontmatter}

\title{Applications of weak convergence for hedging of
game options\thanksref{T1}}
\runtitle{Applications of weak convergence for hedging}

\thankstext{T1}{Supported in part by ISF Grant 130/06.}

\begin{aug}
\author[A]{\fnms{Yan} \snm{Dolinsky}\corref{}\ead[label=e1]{yann1@math.huji.ac.il}}
\runauthor{Y. Dolinsky}
\affiliation{Hebrew University}
\address[A]{Institute of Mathematics\\
Hebrew University of Jerusalem\\
Jerusalem, 91904\\
Israel\\
\printead{e1}} 
\end{aug}

\received{\smonth{3} \syear{2009}}
\revised{\smonth{12} \syear{2009}}

%
\begin{abstract}
In this paper we consider Dynkin's games with payoffs which are
functions of an underlying process. Assuming extended weak convergence
of underlying processes $\{S^{(n)}\}_{n=0}^\infty$ to a limit process
$S$ we prove convergence Dynkin's games values corresponding to
$\{S^{(n)}\}_{n=0}^\infty$ to the Dynkin's game value corresponding to
$S$. We use these results to approximate game options prices with path
dependent payoffs in continuous time models by a sequence of game
options prices in discrete time models which can be calculated by
dynamical programming algorithms. In comparison to previous papers we
work under more general convergence of underlying processes, as well as
weaker conditions on the payoffs.
\end{abstract}

%
\begin{keyword}[class=AMS]
\kwd[Primary ]{91B28}
\kwd[; secondary ]{60F05}
\kwd{91A05}.
\end{keyword}
\begin{keyword}
\kwd{Dynkin games}
\kwd{game options}
\kwd{weak convergence}.
\end{keyword}

\end{frontmatter}

\section{Introduction}\label{sec:1}
Consider a \textit{c\`{a}dl\`{a}g} stochastic process
${\{S_t\}}_{t=0}^T$ ($T<\infty$) which takes on values in
$\mathbb{R}^d_{+}$. For two given functions $f\leq{g}$ let
$X_t=g(t,S)$ and $Y_t=f(t,S)$. Define
$\Gamma(S)=\inf_{\sigma}\sup_{\tau}E(X_{\sigma}\mathbb{I}_{\sigma<\tau
}+Y_{\tau}\mathbb{I}_{\tau\leq\sigma})$
which is the Dynkin game value for the above processes where
$\mathbb{I}_A=1$ if an event $A$ occurs and \mbox{$=$0} if not. The $\inf$
and the $\sup$ are taken over the set of stopping times no bigger
than $T$, with respect to the usual filtration generated by the
process $S$. Our goal is to prove (under some additional
assumptions) that if a sequence of stochastic processes
${\{S^{(n)}_t\}}{}^T_{t=0}$, $n\geq{1}$ converges in law to $S$ then
$\Gamma(S)=\lim_{n\rightarrow\infty}\Gamma(S^{(n)})$.

Although several papers dealt with stability of optimal stopping values under
weak convergence of the underlying processes (see \cite{A2,AK,CKW,CT,L,LP1,LP}
and \cite{MP})
stability of Dynkin's games values under weak convergence of the
underlying processes was not studied before. In his unpublished paper
\cite{A2} Aldous represented the notion of extended weak convergence
and proved the stability of optimal stopping values under extended
weak convergence of the underlying processes. In this paper we extend
these results for Dynkin's
games. The main tools that we use for proving the result is the
Skorohod representation theorem (see \cite{D}) and the theory of
extended weak convergence that Aldous developed in \cite{A2}.

In \cite{CKW} and \cite{L} the authors studied binomial approximations
of American put options in the Black--Scholes (BS) model.
In both of these papers it was proved by using different methods, that
the option price in the BS model is a limit of option prices for an
appropriate sequence of Cox--Ross--Rubinstein (CRR) models. Furthermore
they proved stability of the critical prices and the optimal stopping
times. In
\cite{CT} the authors studied stability of optimal stopping values
and optimal stopping times under convergence in probability of the
underlying processes (though, it seems that the part related to
convergence of optimal stopping values contains some gaps). The
most general result was obtained in \cite{MP} where the authors
considered a general framework and
proved Snell envelopes stability under weak convergence of the
underlying processes. The authors
used the fact that the Snell envelope of a positive process is a
positive supermartingale, and so it is a quasi martingale and the
corresponding tightness theorems (see \cite{MZ}) can be employed.
For Dynkin's games the value process should not, in general, be a
quasi martingale and so the above method is not applicable here, and
so Dynkin's games value process stability under weak convergence
remains an open question, but we are able to prove in this paper the
stability of Dynkin's game values under the weak convergence of the
underlying processes.

One of the motivations to study Dynkin's games values stability under weak
convergence of the underlying processes is applications to game options
approximations. Recall,
that a game contingent claim (GCC) or a game option was defined in
\cite{Ki1} as a contract between the seller and the buyer of the
option such that both have the right to exercise it any time up to a
maturity date (horizon) $T$. If the buyer exercises the contract at
time $t$ then he gets the payment $Y_t$, but if the seller cancels
before the buyer then the latter gets the payment $X_t$. The
difference $\delta_t=X_t-Y_t$ is the penalty which the seller pays
to the buyer for the contract cancellation. Thus the process $S$ can
be considered as a discounted risky asset, and the processes $X\geq
Y$ are considered as the discounted payoff processes. In \cite{Ki1}
it was proved that pricing a GCC in a complete market leads to the
value of a Dynkin game with the payoffs $X,Y$ under the unique
martingale measure, namely if the process $S$ is a martingale then
$\Gamma(S)$ is the option price of the above game option. In
\cite{KK} it was proved that for a general incomplete market, if the
process $S$ is a martingale then $\Gamma(S)$ is an arbitrage-free
price.

Convergence results for Dynkin's games will allow us to approximate
options prices in continuous time markets by a sequence of game
options prices in discrete time markets which are defined on a discrete
probability space. In addition to the theoretical interest such
results have a practical value for calculations of options prices,
since it is well known (see \cite{O}) that for a discrete probability
space Dynkin's games values can be calculated by dynamical
programming algorithm. In this paper we give an example for
approximations of game options with Russian (path dependent) type of
payoffs in the Merton model.

Several papers (see \cite{DK3,Ki2,Ki3}) dealt with
approximations of option prices for game options. These papers used
strong approximation theorems in order to obtain error estimates for
discrete time approximations of game options in the BS model with
Lipschitz conditions on the payoffs. The weak convergence approach
does not allow one to obtain estimates of the error, but it works under
weaker assumptions on the payoffs and can be applied for jump
diffusion models.

The main results of this paper are formulated in the next section where
we also introduce the notation that will be used. In Section \ref{sec3} we
derive auxiliary lemmas that we use. In Section \ref{sec4} we
complete the proof of the main results of the paper. In Section \ref{sec5} we
provide an application for approximations of game options in Merton's
model with path dependent payoffs.

\section{Preliminaries and main results}\label{sec:2}
First we introduce some definitions and notation that will be used
in this paper. Let $d\in\mathbb{N}$. Given a probability space
$(\Omega,\mathcal{F},P)$ consider a \textit{c\`{a}dl\`{a}g}
stochastic process
$S=\{S_t\dvtx\Omega\rightarrow{\mathbb{R}^d}\}_{t=0}^T$. Denote by
$\mathcal{F}^S=\{\mathcal{F}^S_t\}_{t=0}^T$ the usual filtration
of $S$, that is, the smallest right continuous filtration with respect to
which $S$ is adapted, and such that the $\sigma$ algebras contain the
null sets. Let $\mathcal{T}^S_{[0,T]}$ be the set of all stopping
times with respect to $\mathcal{F}^S$ which take values in $[0,T]$.
Denote by $\mathbb{P}^S$ the distribution of $S$ on the canonical
space $\mathbb{D}([0,T];\mathbb{R}^d)$ equipped\vspace*{1pt} with the Skorohod
topology, that is, for any Borel set $A\subset
\mathbb{D}([0,T];\mathbb{R}^d)$, $\mathbb{P}^S(A)=P\{S\in{A}\}$. For
a sequence of stochastic processes
$S^{(n)}\dvtx\Omega_n\rightarrow\mathbb{D}([0,T];\mathbb{R}^d)$ we will
use the notation $S^{(n)}\Rightarrow{S}$ to indicate that the
probability measures $\mathbb{P}^{S^{(n)}}$, $n\geq{1}$ converge
weakly to $\mathbb{P}^S$ (where the space
$\mathbb{D}([0,T];\mathbb{R}^d)$ is equipped with the Skorohod
topology).

Next, let
$f,g\dvtx[0,T]\times\mathbb{D}([0,T];\mathbb{R}^d)\rightarrow\mathbb{R}_{+}$
be two measurable functions such that $f\leq g$. We will assume the
following.
\begin{asm}\label{asm2.1}
For any $t\in[0,T]$ and $x,y\in\mathbb{D}([0,T];\mathbb{R}^d)$
$f(t,x)=f(t,y)$ and $g(t,x)=g(t,y)$ if $x(s)=y(s)$ for any $s\leq{t}$.
The functions
$f(\cdot,x)$, $g(\cdot,x)$ are right-continuous functions with left-hand
limits. Furthermore, let
$\{x_n\}_{n=1}^\infty\subset\mathbb{D}([0,T];\mathbb{R}^d)$ and $\{t_n\}
_{n=1}^\infty\subset{[0,T]}$ such that
$\lim_{n\rightarrow\infty}x_n=x$,
$\lim_{n\rightarrow\infty}t_n=t$ and
$\lim_{n\rightarrow\infty}x_n(t_n)=x(t)$ for some
$x\in\mathbb{D}([0,T];\mathbb{R}^d)$ and $t\in{[0,T]}$. Then
%
\begin{equation}\label{2.1--}
\lim_{n\rightarrow\infty}f(t_n,x_n)=f(t,x)  \quad\mbox{and}\quad
\lim_{n\rightarrow\infty}g(t_n,x_n)=g(t,x).
\end{equation}
\end{asm}

For any \textit{c\`{a}dl\`{a}g} stochastic process
$S=\{S_t\}_{t=0}^T$ set the \textit{c\`{a}dl\`{a}g} adapted
processes
%
\begin{equation}\label{2.1-}
X^{S}_t=g(t,S),\qquad Y^{S}_t=f(t,S)
\end{equation}
and consider the payoff function
%
\begin{equation}\label{2.1}
H^S(t,s)=X^{S}_t\mathbb{I}_{t<s}+Y^{S}_s\mathbb{I}_{s\leq t},\qquad
t,s\leq T.
\end{equation}
Assume that
$\sup_{\tau\in\mathcal{T}^S_{[0,T]}}EX^{S}_{\tau}<\infty$ where $E$
denotes the expectation with respect to the probability measure of
the space on which the process $S$ is defined. Let $\Gamma(S)$ be the
Dynkin's game value of the payoff given by (\ref{2.1}). Namely,
%
\begin{equation}\label{2.2}
\Gamma(S)=\inf_{\sigma\in\mathcal{T}^S_{[0,T]}}\sup_{\tau\in
\mathcal{T}^S_{[0,T]}} E H^S(\sigma,\tau)=\sup_{\tau\in
\mathcal{T}^S_{[0,T]}}\inf_{\sigma\in\mathcal{T}^S_{[0,T]}}E
H^S(\sigma,\tau).
\end{equation}
The second equality follows from Corollary 12 in \cite{LM}.
Furthermore from Lem\-ma~5 in \cite{LM} it follows that
%
\begin{eqnarray}\label{2.3}
\Gamma(S)&=&\inf_{\sigma\in\mathcal{T}^S_{[0,T]}}\sup_{\tau\in
\mathcal{T}^S_{[0,T]}} E J^S(\sigma,\tau)\nonumber\\
&=&\sup_{\tau\in \mathcal{T}^S_{[0,T]}}\inf_{\sigma\in
\mathcal{T}^S_{[0,T]}} EJ^S(\sigma,\tau)\hspace*{120pt}\\
\eqntext{\mbox{where }  J^S(t,s)=\mathbb{I}_{t\wedge s=T}Y^{S}_T+\mathbb{I}_{t\wedge s<T}\mathbb{I}_{t\leq
s}X^{S}_t+Y^{S}_s\mathbb{I}_{s<t},} \\
\eqntext{t,s\leq{T}.}
\end{eqnarray}
First we show that $\Gamma(S)$ depends only on the distribution of
$S$. Consider the probability space
$(\mathbb{D}([0,T];\mathbb{R}^d),P^S)$. Let $\{U_t\}_{t=0}^T$ be
the canonical process of coordinate projection, namely
$U_t\dvtx\mathbb{D}([0,T];\mathbb{R}^d)\rightarrow{\mathbb{R}^d}$ is
given by $U_t(x)=x(t)$ and let $\{\mathcal{G}_t\}_{t=0}^T$ be the
usual filtration which is generated by the above process. Introduce the
set $\Phi$ of all functions
$\phi\dvtx\mathbb{D}([0,T];\mathbb{R}^d)\rightarrow{[0,T]}$ which
satisfy $\{\phi\leq{t}\}\in\mathcal{G}_t$ for any $t\leq{T}$.
Observe that $\sigma\in\mathcal{T}^S_{[0,T]}$ if and only if there
exists a function $\phi\in\Phi$ such that $\sigma=\phi(S)$ a.s. Thus
from (\ref{2.2}) we obtain
%
\begin{equation}\label{2.4}\qquad
\Gamma(S)=\inf_{\phi\in\Phi}\sup_{\psi\in\Phi} E
H^S(\phi(S),\psi(S))= \inf_{\phi\in\Phi}\sup_{\psi\in\Phi}
\mathbb{E}^S H^U(\phi(U),\psi(U)),
\end{equation}
where $\mathbb{E}^S$ is the expectation with respect to the
probability measure $\mathbb{P}^S$. From (\ref{2.4}) it follows that
$\Gamma(S)$ depends only on the distribution of $S$.

In \cite{A2} Aldous introduced the notion of ``extended weak
convergence'' via prediction processes. For the case where the
stochastic processes are considered with respect to their natural
filtration (with the usual assumptions) he proved that extended weak
convergence is equivalent to a more elementary condition which does
not require the use of prediction processes (see \cite{A2}, Proposition
16.15). Following \cite{CT} we will use the above condition as the
definition of extended weak convergence.
\begin{dfn}\label{dfn2.1}
A sequence
$S^{(n)}\dvtx\Omega_n\rightarrow\mathbb{D}([0,T];\mathbb{R}^d)$,
$n\geq{1}$ extended weak converges to a
stochastic process
$S\dvtx\Omega\rightarrow\mathbb{D}([0,T];\mathbb{R}^d)$ if for any $k$
and continuous bounded functions
$\psi_1,\ldots,\psi_k\in{C(\mathbb{D}([0,T];\mathbb{R}^d))}$
%
\begin{equation}\label{2.5}\qquad
\bigl(S^{(n)},Z^{(n,1)},\ldots,Z^{(n,k)}\bigr)\Rightarrow\bigl(S,Z^{(1)},\ldots,Z^{(k)}\bigr)
\qquad \mbox{in }  \mathbb{D}([0,T];\mathbb{R}^{d+k}),
\end{equation}
where for any $t\leq{T}$, $1\leq i\leq{k}$, and $n\in\mathbb{N}$
%
\begin{equation}\label{2.6}\qquad
Z^{(n,i)}_t=E_n\bigl(\psi_i\bigl(S^{(n)}\bigr)|\mathcal{F}^{S^{(n)}}_t\bigr),\qquad
n\in\mathbb{N}  \quad\mbox{and}\quad
Z^{(i)}=E(\psi_i(S)|\mathcal{F}^S_t),
\end{equation}
$E_n$ denotes the expectation with respect to the probability
measure on $\Omega_n$ and $E$ denotes the expectation with respect
to the probability measure on $\Omega$. We will denote extended weak
convergence by $S^{(n)}\Rrightarrow{S}$.
\end{dfn}

Next, we introduce two additional assumptions that we will work
with. Let
$S^{(n)}\dvtx\Omega_n\rightarrow\mathbb{D}([0,T];\mathbb{R}^d)$,
$n\geq{1}$ be a sequence of stochastic processes which satisfies
the following assumptions.
\begin{asm}\label{asm2.2}
The random variables $g(\tau,S^{(n)})$, for $n\geq{1}$ and
$\tau\in\mathcal{T}^{S^{(n)}}_{[0,T]}$ are uniformly integrable.
\end{asm}
\begin{asm}\label{asm2.3}
For any $\varepsilon>0$
\[
\lim_{\delta\downarrow{0}}\lim_{n\rightarrow\infty}\sup_{0<u<\delta}\sup
_{\tau\in\mathcal{F}^{S^{(n)}}_{[0,T]}}
P\bigl(\bigl|S^{(n)}_{(\tau+u)\wedge{T}}-S^{(n)}_{\tau}\bigr|>\varepsilon\bigr)=0.
\]
\end{asm}

The above assumption is called the ``Aldous tightness criterion'' and was
introduced in \cite{A1}. The following theorem is the main result of
the paper.
\begin{theorem}\label{thm2.1}
If $S^{(n)}\Rrightarrow{S}$, then
$\Gamma(S)=\lim_{n\rightarrow\infty}\Gamma(S^{(n)})$.
\end{theorem}

\section{Auxiliary lemmas}\label{sec3}
Let $I\subset{[0,T]}$. For any stochastic process
$S\dvtx\Omega\rightarrow\mathbb{D}([0,T];\mathbb{R}^d)$ denote by
${\Delta}^{S}_{I}$ the set of all stopping times $\tau$ which take
on a finite number of values in $I$ and such that for any $t\in I$,
\mbox{$\{\tau=t\} \in\sigma\{S_u|u\leq{t}\}$.}

\begin{lem}\label{lem3.0+}
Let $I\subset{[0,T]}$ be a dense set which contains the point $T$.
Then
%
\begin{equation}\label{3.0}
\sup_{\tau\in\Delta^S_I}\inf_{\sigma\in\mathcal{F}^S_{[0,T]}}EH^S(\sigma
,\tau)=\Gamma(S)=
\inf_{\sigma\in\Delta^S_I}\sup_{\tau\in\mathcal{F}^S_{[0,T]}}EJ^S(\sigma
,\tau).
\end{equation}
\end{lem}

\begin{pf}
We start with the proof of the first equality. Choose $\varepsilon>0$
and $\tau\in\mathcal{T}^{S}_{[0,T]}$ which satisfies
%
\begin{equation}\label{3.0+}
\inf_{\sigma\in\mathcal{F}^S_{[0,T]}}EH^S(\sigma,\tau)>\Gamma
(S)-\varepsilon.
\end{equation}
For any $n$ let $E_n\subset{I}$ be a finite set which contains $T$
and satisfies $\bigcup_{e\in E_n} (e-\frac{1}{n},e]\supseteq[0,T]$.
Define
%
\begin{equation}\label{3.0++}
\tau_n=\min\biggl\{e\in E_n\Big| e\geq
T\wedge\biggl(\frac{1}{n}+\tau\biggr)\biggr\},\qquad
n\in\mathbb{N}.
\end{equation}
Fix $n$ and $t\in I\setminus{\{T\}}$. Clearly,
\[
\{\tau_n=t\}=\biggl\{\min\biggl\{e\in E_n\Big| e-\frac{1}{n}\geq\tau\biggr\}=t\biggr\}\in
\mathcal{F}^S_{t-{1}/{n}}\subset\sigma\{S_u|u\leq{t}\},
\]
which means that for any $n$, $\tau_n\in{\Delta}^{S}_{I}$. Since
$\tau\leq\tau_n\leq\tau+\frac{2}{n}$, we have $\tau_n\downarrow\tau$.
Set,
$\Phi_n=\sup_{\sigma\in\mathcal{F}^S_{[0,T]}}(H^S(\sigma,\tau
)-H^S(\sigma,\tau_n))^{+}$,
$n\in\mathbb{N}$. Observe that
\[
\limsup_{n\rightarrow\infty}\Phi_n\leq\limsup_{n\rightarrow\infty}\sup
_{\tau\leq
t\leq\tau_n}(Y^S_{\tau}-Y^S_{t})^{+}=0  \qquad\mbox{a.s.}
\]
Since the sequence $\Phi_n$, $n\in\mathbb{N}$ is uniformly
integrable then $\lim_{n\rightarrow\infty}E\Phi_n=0$. From
(\ref{3.0+}) we obtain
\begin{eqnarray*}
&&\Gamma-\sup_{\tau\in\Delta^S_I}\inf_{\sigma\in\mathcal
{F}^S_{[0,T]}}EH^S(\sigma,\tau)\nonumber\\
&&\qquad\leq\varepsilon+\limsup_{n\rightarrow\infty}\sup_{\sigma\in\mathcal{F}^S_{[0,T]}}
E\bigl(H^S(\sigma,\tau)-H^S(\sigma,\tau_n)\bigr)^{+}\nonumber\\
&&\qquad\leq
\limsup_{n\rightarrow\infty}E\Bigl[\sup_{\sigma\in\mathcal
{F}^S_{[0,T]}}\bigl(H^S(\sigma,\tau)-H^S(\sigma,\tau_n)\bigr)^{+}\Bigr]
\\
&&\qquad\leq\varepsilon+\lim_{n\rightarrow\infty}E{\Phi_n}=\varepsilon,
\end{eqnarray*}
and the first equality in (\ref{3.0}) follows. The proof of the
second equality is similar.
\end{pf}

Let $S^{(n)}\Rightarrow{S}$ and assume that the sequence $S^{(n)}$,
$n\in\mathbb{N}$, satisfies the assumptions from Section \ref{sec:2}. The
following two lemmas are a small modification of similar results
that were obtained in \cite{A2}. For the reader's convenience we provide
a self-contained proof for Lemmas \ref{lem3.1} and \ref{lem3.2} which
follows the ideas that were used in~\cite{A2}.
\begin{lem}\label{lem3.1}
Let $E\subset[0,T]$ be a finite set such that any $t\in
E\setminus{\{T\}}$ is a continuity point of the process $S$ a.s.
Then for any $\tau\in{\Delta}^{S}_{E}$ there exists a sequence of
stopping times $\tau_n\in\mathcal{T}^{S^{(n)}}_{[0,T]}$ with values
in $E$ such that $(S^{(n)},\tau_n)\Rightarrow(S,\tau)$ on the space
$\mathbb{D}([0,T];\mathbb{R}^d)\times[0,T]$.
\end{lem}
\begin{pf}
By using the Skorohod representation theorem (see \cite{D}) it follows
that without loss of generality we can assume that there exists a
probability space $(\Omega,\mathcal{F},P)$ on which the process $S$
and the sequence $S^{(n)}$ are defined and $S^{(n)}\rightarrow{S}$
a.s. on $\mathbb{D}([0,T];\mathbb{R}^d)$. In order to prove the
lemma it is sufficient to show that for any $\varepsilon>0$ there
exists $N\in\mathbb{N}$ and a sequence of stopping times
$\tau_n\in\mathcal{T}^{S^{(n)}}_{[0,T]}$ with values in $E$ such
that
%
\begin{equation}\label{3.1}
P \Biggl(\bigcup_{n=N}^{\infty} {\{\tau_n\neq\tau\}}\Biggr)<\varepsilon.
\end{equation}
Choose $\varepsilon>0$. Let $E\setminus{\{T\}}=\{t_1<t_2<\cdots<t_k\}$.
Denote $A_i={\{\tau=t_i\}}\in\sigma\{S_u|u\leq{t_i}\}$, $i\leq{k}$.
Since $t_i$ is a continuity point of the process $S$ there exist
continuous functions
$\phi_i\dvtx\mathbb{D}([0,T];\mathbb{R}^d)\rightarrow{[0,1]}$,
$i\leq{k}$, such that
%
\begin{equation}\label{3.3}
E|\mathbb{I}_{A_i}-\phi_i(S)|<\frac{\varepsilon}{2^{(i+1)}},
\end{equation}
and the function $\phi_i(x)$ depends only on the restriction of $x$
to the interval $[0,t_i]$. For any $n\in\mathbb{N}$ define the
stopping time $\tau_n\in\mathcal{T}^{S_n}_{[0,T]}$ by
%
\begin{equation}\label{3.4}
\tau_n=T\wedge\min\bigl\{t_i|\phi_i\bigl(S^{(n)}\bigr)>\tfrac{1}{2}\bigr\},
\end{equation}
where $\min\{t_i|\phi_i(S^{(n)})>\frac{1}{2}\}=\infty$ if for any
$i$, $\phi_i(S^{(n)})\leq\frac{1}{2}$. Observe that
$\phi^{(n)}_i(S^{(n)})$ is a $\mathcal{F}^{S_n}_{t_i}$ measurable
random variable, thus $\tau_n$ is indeed a stopping time with
respect to the filtration $\mathcal{F}^{S^{(n)}}$. Set
\[
C_i=\bigl(A_i\cap\{\phi_i(S)>1/2\}\bigr)\cup\bigl(A^{c}_i\cap
\{\phi_i(S)<1/2\}\bigr),\qquad  i\leq{k},
\]
and
\[
C=\bigcap_{i=1}^k
C_i.
\]
Since for any $i$, $\phi_i$ is a continuous function we
obtain that for any $\omega\in{C}$ there exists $N(\omega)$ such
that for any $n\geq{N(\omega)}$, $\tau(\omega)=\tau_n(\omega)$. From
(\ref{3.3}) and the Markov inequality we obtain
%
\begin{eqnarray}\label{3.5}
P(C)&=&1-P(\Omega\setminus C)\geq1-\sum_{i=1}^k
P\biggl(|\phi_i(S)-\mathbb{I}_{A_i}|\geq{\frac{1}{2}}\biggr)\nonumber\\[-8pt]\\[-8pt]
&\geq& 1- \sum_{i=1}^k \frac{\varepsilon}{2^{i}}>1 -\varepsilon.\nonumber
\end{eqnarray}
Set $E_n=\bigcap_{m=n}^{\infty} {\{\tau_m=\tau\}}$,
$n\in\mathbb{N}$. Observe that the sequence $E_n$, $n\geq{1}$ is an
increasing sequence of events and $\bigcup_{n=1}^{\infty}
E_n\supset{C}$. From (\ref{3.5}) it follows that there exists
$N\in\mathbb{N}$ such that $P(E_N)>1-\varepsilon$ and (\ref{3.1})
follows.
\end{pf}

\begin{lem}\label{lem3.2}
Assume that $\tau_n\in\mathcal{T}^{S^{(n)}}_{[0,T]}$, $n\geq{1}$ is
a sequence of stopping times which satisfies
$(S^{(n)},\tau_n)\Rightarrow(S,\nu)$ on the space
$\mathbb{D}([0,T];\mathbb{R}^d)\times[0,T]$ for some random
variable $\nu$. Then
%
\begin{equation}\label{3.5+}
\bigl(S^{(n)},S^{(n)}_{\tau_n},\tau_n\bigr)\Rightarrow(S,S_{\nu},\nu)
\end{equation}
on the space $\mathbb{D}([0,T];\mathbb{R}^d)\times
\mathbb{R}^d\times[0,T]$. In addition, if $S^{(n)}\Rrightarrow{S}$,
then for any $t$, $\{\nu\leq{t}\}$ and $\mathcal{F}^S_T$ are
conditionally independent given $\mathcal{F}^S_t$, and so for any
uniformly integrable \textit{c\`{a}dl\`{a}g} stochastic process
$\{V_t\}_{t=0}^T$ adapted to the filtration
$\mathcal{F}^S_{[0,T]}$
%
\begin{equation}\label{3.6}
\inf_{\tau\in\mathcal{T}^S_{[0,T]}}EV_{\tau}\leq
EV_{\nu}\leq\sup_{\tau\in\mathcal{T}^S_{[0,T]}}EV_{\tau}.
\end{equation}
\end{lem}
\begin{pf}
By using the Skorohod representation theorem it follows that without
loss of generality we can assume that there exists a probability
space $(\Omega,\mathcal{F},P)$ on which the process $S,\nu$ and the
sequence $S^{(n)},\tau_n$ are defined and
$(S^{(n)},\tau_n)\rightarrow(S,\nu)$ a.s. on
$\mathbb{D}([0,T];\mathbb{R}^d)\times[0,T]$. Thus in order to prove
(\ref{3.5+}) it is sufficient to show that
$S^{(n)}_{\tau_n}\rightarrow{S_{\nu}}$ in probability. Choose
$\varepsilon>0$. The process $Z_t:=S_{(\nu+t)\wedge{T}}$, $t\geq{0}$ is
a \textit{c\`{a}dl\`{a}g} process. It is well known (see, e.g.,
\cite{B}, Chapter 3) that for a \textit{c\`{a}dl\`{a}g}
process the set of points for which the process is not continuous
(with positive probability) is at most countable, thus there exists
a sequence $u_n\downarrow{0}$ such that for any $n$ the process $Z$
is continuous at $u_n$, which means that for any $\omega\in\Omega$,
$\nu(\omega)+u_n$ is a continuity point of the function $S(\omega)$
provided that $\nu(\omega)+u_n<T$. Since the map
$(f,t)\rightarrow{f(t)}$ from
$\mathbb{D}([0,T];\mathbb{R}^d)\times{[0,T]}$ to $\mathbb{R}^d$ is
continuous at $(f_0,t_0)$ if $t_0$ is a continuity point of $f_0$
(see \cite{B}, Chapter 3) we obtain that for any
$\omega\in{E_1}:=\{\nu<T\}$,
$\lim_{n\rightarrow\infty}S^{(n)}_{(\tau_n+u_n)\wedge{T}}=S_{\nu}$.
This together with Assumption \ref{asm2.3} gives
%
\begin{eqnarray}\label{3.10}
&&\lim_{n\rightarrow\infty}P\bigl(E_1\cap\bigl\{\bigl|S_{\nu}-S^{(n)}_{\tau
_n}\bigr|>2\varepsilon\bigr\}\bigr)\nonumber\\
&&\qquad
\leq\lim_{n\rightarrow\infty}P\bigl(E_1\cap\bigl\{\bigl|S_{\nu}-S^{(n)}_{(\tau_n+u_n)\wedge
{T}}\bigr|>\varepsilon\bigr\}\bigr)\\
&&\qquad\quad{}+
\lim_{n\rightarrow\infty}P\bigl\{\bigl|S^{(n)}_{\tau_n}-S^{(n)}_{(\tau
_n+u_n)\wedge{T}}\bigr|>\varepsilon\bigr\}=0.\nonumber
\end{eqnarray}
Next, we deal with the event $E_2:=\{\nu=T\}$. For any $\delta>0$
and $n\in\mathbb{N}$ set
$\tau^{(\delta)}_n=\mathbb{I}_{\tau_n<T-\delta}\tau_n+\mathbb{I}_{\tau
_n\geq
T-\delta}T\in\mathcal{T}^{S^{(n)}}_{[0,T]}$. Observe that for any
$\omega\in E_2$ there exists $N(\omega)\in\mathbb{N}$ such that for
any $n>N(\omega)$, $\tau^{(\delta)}_n=T$. Since the map
$f\rightarrow{f(T)}$ from $\mathbb{D}([0,T];\mathbb{R}^d)$ to
$\mathbb{R}^d$ is continuous we obtain that for any $\delta>0$ and
$\omega\in E_2$,
$\lim_{n\rightarrow\infty}S^{(n)}_{\tau^{(\delta)}_n}=S_T=S_{\nu}$.
Thus from Assumption \ref{asm2.3} we obtain
%
\begin{eqnarray}\label{3.12}
&&\lim_{n\rightarrow\infty}P\bigl(B\cap\bigl\{\bigl|S_{\nu}-S^{(n)}_{\tau_n}\bigr|>2\varepsilon
\bigr\}\bigr)\nonumber\\
&&\qquad \leq\limsup_{\delta\downarrow{0}}\limsup_{n\rightarrow\infty}
P\bigl(B\cap\bigl\{\bigl|S_{\nu}-S^{(n)}_{\tau^{(\delta)}_n}\bigr|>\varepsilon\bigr\}\bigr)\\
&&\qquad\quad{}+
\limsup_{\delta\downarrow{0}}\limsup_{n\rightarrow\infty}P\bigl\{
\bigl|S^{(n)}_{\tau_n}-S^{(n)}_{\tau^{(\delta)}_n}\bigr|>\varepsilon\bigr\}=0.
\nonumber
\end{eqnarray}
From (\ref{3.10}) and (\ref{3.12}) it follows that
$\lim_{n\rightarrow\infty}P\{|S_{\nu}-S^{(n)}_{\tau_n}|>2\varepsilon\}=0$
and (\ref{3.5+}) follows. Next, let $S^{(n)}\Rrightarrow{S}$ and
assume without loss of generality that $(\Omega,\mathcal{F},P)$ is
large enough such that there exists a random variable $H$
distributed uniformly on the interval $[0,1]$ and independent of
$\mathcal{F}^S_T$. First we show that for any $t<T$,
$\{\nu\leq{t}\}$ and $\mathcal{F}^S_T$ are conditionally independent
given $\mathcal{F}^S_t$, that is,
%
\begin{equation}\label{3.13}
E(\mathbb{I}_{\nu\leq{t}}|\mathcal{F}^S_t)=E(\mathbb{I}_{\nu\leq
{t}}|\mathcal{F}^S_T).
\end{equation}
Fix $t<T$, a $\psi\in{C(\mathbb{D}([0,T];\mathbb{R}^d))}$ and
$\phi\in{C[0,T]}$. Define the (\textit{c\`{a}dl\`{a}g}) stochastic
processes $Z_{u}=E(\psi(S)|\mathcal{F}^S_u)$ and
$Z^{(n)}_{u}=E(\psi(S^{(n)})|\mathcal{F}^{S^{(n)}}_u)$. Let\vadjust{\goodbreak}
$u_n\downarrow{t}$ be a sequence such that for any $n$ the process
$Z$ is continuous at $u_n$. Clearly,
%
\begin{equation}\label{3.15}
E\bigl[\phi(\tau_m\wedge
u_n)\bigl(\psi\bigl(S^{(m)}\bigr)-Z^{(m)}_{u_n}\bigr)\bigr]=0\qquad
\forall{n,m}\in\mathbb{N}.
\end{equation}
Since $S^{(m)}\Rrightarrow{S}$ we obtain that for any $n$,
$Z^{(m)}_{u_n}\Rightarrow Z_{u_n}$. Fix $n$. The sequence
$(S^{(m)},Z^{(m)}_{u_n},\tau_m)$, $m\geq{1}$ is tight and so from
Prohorov's theorem (see \cite{B}) it follows that there exists a
subsequence $(S^{(m_k)},Z^{(m_k)}_{u_n},\tau_{m_k})$ which converges
in law to $(S,Z_{u_n},\nu)$. This together with (\ref{3.15}) gives
$E[\phi(\nu\wedge u_n) (\psi(S)-Z(u_n))]=0$. The function
$\psi\in{C(\mathbb{D}([0,T];\mathbb{R}^d))}$ is arbitrary, and so
from density arguments it follows that for any
$B\in\mathcal{F}^S_{[0,T]}$ and $n\in\mathbb{N}$, $E[\phi(\nu\wedge
u_n) (\mathbb{I}_B- E(\mathbb{I}_B|\mathcal{F}^{S}_{u_n}))]=0$.
Since $\{\nu\leq{t}\}=\{\nu\wedge u_n\leq t\}$ and $\phi$ is
arbitrary then by using density arguments it follows that
$E[\mathbb{I}_{\nu\leq{t}}
(\mathbb{I}_B-E(\mathbb{I}_B|\mathcal{F}^{S}_{u_n}))]=0$, and by
letting $n\rightarrow\infty$ we obtain that for any
$B\in\mathcal{F}^S_{[0,T]}$, $E[\mathbb{I}_{\nu\leq{t}}
(\mathbb{I}_B-E(\mathbb{I}_B|\mathcal{F}^{S}_{t}))]=0$. Thus for any
$B\in\mathcal{F}^S_{[0,T]}$,
\begin{eqnarray*}
E[(\mathbb{I}_{\nu\leq{t}}|\mathcal{F}^S_T)\mathbb{I}_B]&=&E(\mathbb
{I}_{\nu\leq{t}}\mathbb{I}_B)=
E[\mathbb{I}_{\nu\leq{t}}E(\mathbb{I}_B|\mathcal{F}^S_t)]\\
&=& E[E(\mathbb{I}_{\nu\leq{t}}|\mathcal{F}^S_t)E(\mathbb{I}_B|\mathcal
{F}^{S}_t)]\\
&=&
E[E(\mathbb{I}_{\nu\leq{t}}|\mathcal{F}^S_t)\mathbb{I}_B]
\end{eqnarray*}
and (\ref{3.13}) follows. Next, define the stochastic process
$Q_t=E(\mathbb{I}_{\nu\leq{t}}|\mathcal{F}^S_T)$, $t\leq{T}$.
Clearly $Q$ is a positive increasing (adapted to $\mathcal{F}^S$)
\textit{c\`{a}dl\`{a}g} process and $Q_T=1$ a.s. Set
$\sigma=\inf\{t|Q_t\geq{H}\}$. From the definition of $H$ we obtain
\[
E(\mathbb{I}_{\sigma\leq{t}}|\mathcal{F}^S_T)=E(\mathbb{I}_{Q_t\geq
{H}|\mathcal{F}^S_t})=Q_t=
(\mathbb{I}_{\nu\leq{t}}|\mathcal{F}^S_T),
\]
and since $V$ is $\mathcal{F}^S_T$ measurable, then
$EV_{\sigma}=EV_{\nu}$. Finally, for any $0\leq u\leq1$ define
$\sigma_u=\inf\{t|Q_t\geq{u}\}\in\mathcal{T}^S_{[0,T]}$. Since $H$
is independent of $V$ and $Q$, then $EV_{\nu}=EV_{\sigma}=\int_{0}^1
(EV_{\sigma_u})\,du$ and (\ref{3.6}) follows.
\end{pf}

\section{Proof of main results}\label{sec4}

In this section we complete the proof of Theorem
\ref{thm2.1}. Denote $\Gamma=\Gamma(S)$ and
$\Gamma_n=\Gamma(S^{(n)})$, $n\geq{1}$. First we prove that
$\Gamma\leq\lim_{n\rightarrow\infty}\Gamma_n$. Here and in the
sequel, for the sake of simplicity we will assume that indices have
been renamed so that the whole sequence converges. Choose
$\varepsilon>0$. Denote by $I\subset[0,T]$ the union of the point
$\{T\}$ together with all continuity points of the process $S$. From
Lemma \ref{lem3.0+} it follows that there exists
$\tau\in\Delta^S_{I}$ such that
%
\begin{equation}\label{4.0+}
\Gamma(S)<\varepsilon+\inf_{\sigma\in\mathcal{F}^S_{[0,T]}}EH^S(\sigma,\tau).
\end{equation}
From Lemma \ref{lem3.1} we can choose a sequence of stopping times
$\tau_{n}\in\mathcal{T}^{S^{(n)}}_{[0,T]}$, $n\geq{1}$, such that
$(S^{(n)},\tau_n)\Rightarrow(S,\tau)$ on the space
$\mathbb{D}([0,T];\mathbb{R}^d)\times[0,T]$. From
Lem\-ma~\ref{lem3.2} we obtain
%
\begin{equation}\label{4.1-}
\bigl(S^{(n)},S^{(n)}_{\tau_n},\tau_n\bigr)\Rightarrow(S,S_{\tau},\tau)
\end{equation}
on the space $\mathbb{D}([0,T];\mathbb{R}^d)\times
\mathbb{R}^{d}\times[0,T]$. From (\ref{2.2}) it follows that for
any $n\in\mathbb{N}$ there exists a stopping time
$\sigma_n\in\mathcal{T}^{S^{(n)}}_{[0,T]}$ such that
%
\begin{equation}\label{4.1}
\Gamma_n> E_n H^{S^{(n)}}(\sigma_n,\tau_n)-\varepsilon.
\end{equation}
The sequence $(S^{(n)},\sigma_{n})$ is tight in
$\mathbb{D}([0,T];\mathbb{R}^d)\times{[0,T]}$ and so
$(S^{(n)},\sigma_{n})\Rightarrow{(S,\nu)}$ for some random variable
$\nu\leq{T}$. From Lemma \ref{lem3.2}
%
\begin{equation}\label{4.1+}
\bigl(S^{(n)},S^{(n)}_{\sigma_{n}},\sigma_{n}\bigr)\Rightarrow(S,S_{\nu},\nu)
\end{equation}
on the space $\mathbb{D}([0,T];\mathbb{R}^d)\times
\mathbb{R}^{d}\times[0,T]$. From (\ref{4.1-}) and (\ref{4.1+}) it
follows that the sequence
$(S^{(n)},S^{(n)}_{\tau_{n}},S^{(n)}_{\sigma_{n}},\tau_{n},\sigma_{n})$,
$n\geq{1}$, is tight on the space
$\mathbb{D}([0,T];\mathbb{R}^d)\times\mathbb{R}^{2d}\times
{[0,T]}^2$. Thus
$(S^{(n)},S^{(n)}_{\tau_{n}},S^{(n)}_{\sigma_{n}},\tau_{n},\sigma
_{n})\Rightarrow
(S,S_{\tau},S_{\nu},\tau,\nu)$. By using the Skorohod representation
theorem it follows that without loss of generality we can assume
that there exists a probability space $(\Omega,\mathcal{F},P)$ on
which $(S^{(n)},S^{(n)}_{\tau_{n}}$, $S^{(n)}_{\sigma_{n}},
\tau_{n},\sigma_{n})\rightarrow(S,S_{\tau},S_{\nu},\tau,\nu)$ a.s.
on the space $\mathbb{D}([0,T];\mathbb{R}^d)\times
\mathbb{R}^{2d}\times{[0,T]}^2$. This together with Assumption
\ref{asm2.1} gives
%
\begin{equation}\label{4.3}
H^S(\nu,\tau)\leq\liminf_{n\rightarrow\infty}H^{S^{(n)}}(\sigma_{n},\tau_{n}).
\end{equation}
From (\ref{4.1}) and (\ref{4.3})
%
\begin{equation}\label{4.4-}
EH^S(\nu,\tau)\leq
\liminf_{n\rightarrow\infty}EH^{S^{(n)}}(\sigma_{n},\tau_{n})\leq
\lim_{n\rightarrow\infty}\Gamma_n+\varepsilon.
\end{equation}
By applying Lemma \ref{lem3.2} for the process $Q_t:=H^S(t,\tau)$ it
follows
%
\begin{equation}\label{4.4}
\inf_{\sigma\in\mathcal{F}^S_{[0,T]}}EH^S(\sigma,\tau)\leq
EH^S(\nu,\tau)\leq\lim_{n\rightarrow\infty}\Gamma_n+\varepsilon.
\end{equation}
From (\ref{4.0+}) and (\ref{4.4}) we obtain
$\Gamma\leq\lim_{n\rightarrow\infty}\Gamma_n$. In order to complete
the proof we prove that
$\Gamma\geq\lim_{n\rightarrow\infty}\Gamma_n$. Choose
$\varepsilon\geq{0}$. From Lemma \ref{lem3.0+} there exists a stopping
time $\sigma\in\Delta^S_{I}$ which takes values on a finite set $E$
and satisfies
%
\begin{equation}\label{4.4+}
\Gamma(S)>\sup_{\tau\in\mathcal{F}^S_{[0,T]}}EJ^S(\sigma,\tau)-\varepsilon.
\end{equation}
From Lemma \ref{lem3.1} and \ref{lem3.2} it follows that we can
choose a sequence of stopping times
$\sigma_{n}\in\mathcal{T}^{S^{(n)}}_{[0,T]}$, $n\geq{1}$, with values
in $E$ such that
%
\begin{equation}
\bigl(S^{(n)},S^{(n)}_{\sigma_n},\sigma_n\bigr)\Rightarrow
(S,S_{\sigma},\sigma)
\end{equation}
in law on the space
$\mathbb{D}([0,T];\mathbb{R}^d)\times\mathbb{R}^d\times[0,T]$. From
(\ref{2.3}) it follows that for any $n\in\mathbb{N}$ there exists a
stopping time $\tau_n\in\mathcal{T}^{S^{(n)}}_{[0,T]}$ such that
%
\begin{equation}\label{4.7-}
\Gamma_n<\inf_{\sigma_n\in\mathcal{T}^{S^{(n)}}_{[0,T]}}E_nJ
^{S^{(n)}}(\sigma_n,\tau_n)+\varepsilon.
\end{equation}
The sequence $(S^{(n)},\tau_n)$ is tight, and thus from Lemma
\ref{lem3.2},
%
\begin{equation}\label{4.7}
\bigl(S^{(n)},S^{(n)}_{\tau_n},\tau_n\bigr)\Rightarrow(S,S_{\nu},\nu)
\end{equation}
on the space $\mathbb{D}([0,T];\mathbb{R}^d)\times
\mathbb{R}^{d}\times[0,T]$, for some $\nu\leq{T}$. As before, by
using the Skorohod representation theorem it follows that without loss
of generality we can assume that there exists a probability space
$(\Omega,\mathcal{F},P)$ on which
$(S^{(n)},S^{(n)}_{\tau_{n}},S^{(n)}_{\sigma_{n}},\tau_{n},\sigma
_{n})\rightarrow
(S,S_{\nu},S_{\sigma},\nu,\sigma)$ a.s. on the space
$\mathbb{D}([0,T];\mathbb{R}^d)\times\mathbb{R}^{2d}\times
{[0,T]}^2$. Observe that if $\sigma=T$, then $\sigma_n=T$ for
sufficiently large $n$. Thus from Assumption \ref{asm2.1}
\[
J^S(\sigma,\nu)\geq\limsup_{n\rightarrow\infty}J^{S^{(n)}}(\sigma_n,\tau_n).
\]
This together with (\ref{4.7-}) and Assumption \ref{asm2.2} gives
%
\begin{equation}\label{4.10}
EJ^S(\sigma,\nu)\geq
\limsup_{n\rightarrow\infty}EJ^{S^{(n)}}(\sigma_n,\tau_n)\geq
\lim_{n\rightarrow\infty}\Gamma_n-\varepsilon.
\end{equation}
By applying Lemma \ref{lem3.2} for the process $Q_t:=J^S(\sigma,t)$
we obtain
%
\begin{equation}\label{4.11}
\sup_{\tau\in\mathcal{T}^S_{[0,T]}}EJ^{S}(\sigma,\tau)\geq
EJ^{S}(\sigma,\nu).
\end{equation}
From (\ref{4.4+}), (\ref{4.10}) and (\ref{4.11}),
$\Gamma\geq\lim_{n\rightarrow\infty}\Gamma_n$.

\section{Applications to game options}\label{sec5}

In this section we give an example for an application of Theorem
\ref{thm2.1}. We will consider discrete time approximations of game
options prices in the Merton (one-dimensional) model. Approximation
of American options in the Merton model were considered in \cite{M1}.
Let $(\Omega,\mathcal{F},P)$ be a probability space together with a
standard Brownian motion $\{W_t\}_{t=0}^T$, a Poisson process
$\{N_t\}_{t=0}^T$ with intensity $\lambda$ and independent of $W$
and a sequence of i.i.d. random variables ${\{U_i\}}_{i=1}^\infty$
with values in $(-1,\infty)$, independent of $N$ and $W$. We
assume that $EU_1<\infty$. A Merton model with horizon $T<\infty$
consists of a savings account with an interest rate $r>0$, and of a
risky asset (stock). Assume that the discounted stock price
$\{S_t\}_{t=0}^T$ [i.e., a ratio of the original stock price and
$\exp(rt)$] is given by
%
\begin{eqnarray}\label{5.1}
S_t=S_0\exp\Biggl(\biggl(\mu-\frac{{\sigma}^2}{2}\biggr)t+\sigma W_t+\sum_{j=1}^{N_t}
\ln(1+U_j)\Biggr),\nonumber\\[-8pt]\\[-8pt]
\eqntext{S_0,\sigma>0,  \mu=-\lambda EU_1.}
\end{eqnarray}
The equality $\mu=-\lambda EU_1$ guarantees that $S$ is a martingale
with respect to $P$ and the usual filtration $\mathcal{F}^S$.
Introduce\vadjust{\goodbreak} a game option with Russian payoff functions. Namely the
discounted payoffs are $Y^S_t=f(t,S)$ and $X^S_t=g(t,S)$ where
$f,g\dvtx[0,T]\times\mathbb{D}([0,T];\mathbb{R})\rightarrow\mathbb{R}_{+}$
are given by
%
\begin{eqnarray}\label{5.2}
f(t,x)&=&e^{-rt}\max\Bigl(M,\sup_{0\leq u\leq t}e^{ru}x_u \Bigr)
\quad\mbox{and}\nonumber\\[-8pt]\\[-8pt]
g(t,x)&=&f(t,x)+\delta x_t,\qquad  \delta,r,M>0.\nonumber
\end{eqnarray}
From \cite{KK} it follows that
%
\begin{equation}\label{5.3}
V:=\Gamma(S)=\inf_{\sigma\in\mathcal{T}^S_{[0,T]}}\sup_{\tau\in
\mathcal{T}^S_{[0,T]}} E H^S(\sigma,\tau)
\end{equation}
is an arbitrage-free price (recall that the Merton model is
incomplete). Following \cite{M1} we construct a sequence of discrete
time approximations. For any $n\in\mathbb{N}$ let
$(\Omega_n,\mathcal{F}_n,P_n)$ be a probability space together with
three independent sequences of i.i.d. random variables
$\{\xi^{(n)}_k\}_{k=1}^n$, $\{\rho^{(n)}_k\}_{k=1}^n$ and
$\{u^{(n)}_k\}_{k=1}^n$. The first one is a sequence of
Bernoulli random variables such that
$P_n\{\xi^{(n)}_k=1\}=\frac{\lambda T}{n}$, the second sequence
satisfies $P_n\{\rho^{(n)}_k=1\}=1-P_n\{\rho^{(n)}_k=-1\}=
\frac{{n}/({n+\mu\lambda
T})-\exp(-\sigma\sqrt{T/n})}{\exp(\sigma\sqrt{T/n})-\exp
(-\sigma\sqrt{T/n})}$
(we assume that $n$ is sufficiently large such that above term is
positive) and the third sequence given by $u^{(n)}_k\sim
\ln(1+U_1)$. For any $0\leq k\leq{n}$ and $kT/n\leq t<(k+1)T/n$ set
%
\begin{eqnarray}\label{5.4}
W^{(n)}_t&=&\sqrt\frac{T}{n}\sum_{i=1}^k \rho^{(n)}_i,\qquad
N^{(n)}_t=\sum_{i=1}^k
\xi^{(n)}_i  \quad\mbox{and} \nonumber\\[-8pt]\\[-8pt]
S^{(n)}_t&=&S_0\exp\Biggl(\sigma W^{(n)}_t+\sum_{i=1}^{N^{(n)}_t}
u^{(n)}_i\Biggr). \nonumber
\end{eqnarray}
The $n$-step discrete time market is active at times
$\{0,\frac{T}{n},\frac{2T}{n},\ldots,T\}$ and consists of a savings
account with an interest rate $r>0$, and of a risky asset whose
discounted stock price $S^{(n)}$ is given by (\ref{5.4}). Consider a
game option with the discounted payoffs $Y^{S^{(n)}}_t=f(t,S^{(n)})$ and
$X^{S^{(n)}}_t=g(t,S^{(n)})$. Let $\Delta_n$ be the set of all
stopping times with respect to the filtration
$\mathcal{F}^{S^{(n)}}$ with values in the set
$\{0,\frac{T}{n},\frac{2T}{n},\ldots,T\}$. Since the process
$\{S^{(n)}_{kT/n}\}_{k=0}^n$ is a martingale under $P_n$ it
follows that $V_n$ which is given by
%
\begin{equation}\label{5.6}
V_n=\inf_{\sigma\in\Delta_n}\sup_{\tau\in\Delta_n} E_n
H^{S^{(n)}}(\sigma,\tau)
\end{equation}
is an arbitrage-free price. Next, we describe a dynamical
programming algorithm which allows us to calculate $V_n$.\vadjust{\goodbreak} For $0\leq
k\leq{n}$ define the functions
$\psi^{(n)}_k,\phi^{(n)}_k\dvtx\mathbb{R}^k\times\{-1,1\}^k\times
\{0,1\}^k\rightarrow\mathbb{R}_{+}$ by
%
\begin{eqnarray}\label{5.7}
&&\psi^{(n)}_k(x_1,\ldots,x_k,y_1,\ldots,y_k,z_1,\ldots,z_k)\hspace*{-18pt}\nonumber\\
&&\qquad=\exp\biggl(-\frac
{rkT}{n}\biggr)\nonumber\\
&&\qquad\quad{}\times\max \Biggl(M,S_0\max_{0\leq{i}\leq{k}}\exp\Biggl(\frac{r ((i+1)\wedge{k})
T}{n}\hspace*{-18pt}\nonumber\\
&&\qquad\quad\hspace*{161pt}\hspace*{-48pt}{}+\sigma\sqrt\frac{T}{n}\sum_{j=1}^i y_j+\sum_{j=1}^{m_i} x_j
\Biggr)\Biggr) \quad\mbox{and} \hspace*{-18pt}\\
&&\phi^{(n)}_k(x_1,\ldots,x_k,y_1,\ldots,y_k,z_1,\ldots,z_k)\hspace*{-18pt}
\nonumber\\
&&\qquad=\psi^{(n)}_k(x_1,\ldots,x_k,y_1,\ldots,y_k,z_1,\ldots,z_k)\hspace*{-18pt}\nonumber\\
&&\qquad\quad{}+\delta
S_0\exp\Biggl(\sigma\sqrt\frac{T}{n}\sum_{i=1}^k y_i+\sum_{i=1}^{m_k}
x_i\Biggr),\hspace*{-18pt}\nonumber
\end{eqnarray}
where $m_i=m_i(z_1,\ldots,z_k)=\sum_{q=1}^i z_q$. Since the process
$S^{(n)}$ is constant on intervals of the form $[iT/n,(i+1)T/n)$,
for any $0\leq{k}\leq{n}$
%
\begin{eqnarray}\label{5.8}
\psi^{(n)}_k\bigl(u^{(n)}_1,\ldots,u^{(n)}_k,\rho^{(n)}_1,\ldots,\rho^{(n)}_k,\xi
^{(n)}_1,\ldots,\xi^{(n)}_k\bigr)&=&Y^{S^{(n)}}_{{kT}/{n}}\quad\mbox{and}\nonumber\\[-8pt]\\[-8pt]
\phi^{(n)}_k\bigl(u^{(n)}_1,\ldots,u^{(n)}_k,\rho^{(n)}_1,\ldots,\rho^{(n)}_k,\xi
^{(n)}_1,\ldots,\xi^{(n)}_k\bigr)&=&X^{S^{(n)}}_{{kT}/{n}}.\nonumber
\end{eqnarray}
Finally, define a sequence $\{J^{(n)}_k\}_{k=0}^n$ of functions
$J^{(n)}_k\dvtx\mathbb{R}^k\times\{-1,1\}^k\times
\{0,1\}^k\rightarrow\mathbb{R}_{+}$ by the following backward
recursion:
%
\begin{eqnarray}\label{5.9}
&&\hspace*{1.4pt}J^{(n)}_n=\psi^{(n)}_n   \quad\mbox{and}\nonumber\\
&&J^{(n)}_k(x_1,\ldots,x_k,y_1,\ldots,y_k,z_1,\ldots,z_k)\nonumber\\
&&\qquad=\min \bigl(\phi
^{(n)}_k(x_1,\ldots,x_k,y_1,\ldots,y_k,z_1,\ldots,z_k),\nonumber\\
&&\qquad\hspace*{32.9pt}\max \bigl(\psi^{(n)}_k(x_1,\ldots,x_k,y_1,\ldots,y_k,z_1,\ldots,z_k),\\
&&\qquad\hspace*{56.5pt}E_n\bigl[J^{(n)}_{k+1}\bigl(x_1,\ldots,x_k,u^{(n)}_1,\nonumber\\
&&\qquad\hspace*{99pt}y_1,\ldots,y_k,\rho^{(n)}_1,z_1,\ldots,z_k,\xi^{(n)}_1\bigr)\bigr]\bigr)
\bigr)\nonumber\\
\eqntext{\mbox{for } k=n-1,n-2,\ldots,0.}
\end{eqnarray}
From (\ref{5.8}) and by using the dynamical programming algorithm
that was obtained in \cite{O} for general Dynkin's games in discrete
time, it follows that
%
\begin{equation}\label{5.10}
V_n=J^{(n)}_0.
\end{equation}
The following result says that the option price in the continuous
time Merton model can be approximated by the sequence
${\{V_n\}}_{n=1}^\infty$ (which can be calculated explicitly as
shown above).
\begin{theorem}\label{thm5.1}
$V=\lim_{n\rightarrow\infty}V_n$.
\end{theorem}
\begin{pf}
First we prove that $V=\lim_{n\rightarrow\infty}\Gamma(S^{(n)})$.
Let us check that the conditions of Theorem \ref{thm2.1} are
satisfied. It can be easily checked that the functions $f,g$
satisfy Assumption \ref{asm2.1} and that for any $k>1$, $\sup_{n\geq
1} E_n(S^{(n)}_T)^k
<\infty$. Thus from Doob's inequality we obtain
that $\sup_{n\geq1}$
$E_n[(\sup_{0\leq{t}\leq{T}}S^{(n)}_t)^k]
<\infty$ and Assumption
\ref{asm2.2} follows. It is well known that
$(W^{(n)},N^{(n)})\Rightarrow(W,N)$ on the space
$\mathbb{D}([0,T];\mathbb{R}^{2})$, and so by using the Skorohod
representation theorem we can build a probability space on which
$(W^{(n)},N^{(n)})\rightarrow(W,N)$ a.s. and on which there exists a
sequence of i.i.d. random variables $U_1,\ldots,U_n,\ldots$ which is
independent of $\{W^{(n)},N^{(n)}\}_{n=1}^\infty,W,N$. Thus $\sigma
W^{(n)}+\sum_{i=1}^{N^{(n)}} U_i\rightarrow\sigma W+\sum_{i=1}^{N}
U_i$ a.s. on the space $\mathbb{D}([0,T];\mathbb{R})$ thus $\ln
S^{(n)}\Rightarrow\ln S$. For any $n$ the process $\ln S^{(n)}$ has
independent increments and the process $\ln S$ is continuous in
probability with independent increments. From Corollaries 1 and 2 in
\cite{JS} we obtain that $\ln S^{(n)}\Rrightarrow\ln S$ and that
$\ln S^{(n)}$, $n\geq{1}$, satisfies Assumption \ref{asm2.3}, which
means that $S^{(n)}\Rrightarrow S$ and $S^{(n)}$, $n\geq{1}$,
satisfies Assumption \ref{asm2.3}. We conclude that the conditions
of Theorem \ref{thm2.1} are satisfied, and the equality
$V=\lim_{n\rightarrow\infty}\Gamma(S^{(n)})$ follows. In order to
complete the proof it remains to show that
%
\begin{equation}\label{5.11}
\lim_{n\rightarrow\infty}\bigl|V_n-\Gamma\bigl(S^{(n)}\bigr)\bigr|=0.
\end{equation}
For any $n$ define the maps
$\Phi_n,\Psi_n\dvtx\mathcal{T}^{S^{(n)}}_{[0,T]}\rightarrow\Delta_n$
%
\begin{equation}\label{5.12}\hspace*{28pt}
\Phi_n(\sigma)=\frac{T}{n}\max\{k|kT/n\leq\sigma\}
\quad\mbox{and}\quad
\Psi_n(\sigma)=\frac{T}{n}\min\{k|kT/n\geq\sigma\}.
\end{equation}
From (\ref{5.6})
\[
\inf_{\sigma\in
\Delta_n}\sup_{\tau\in\mathcal{T}^{S^{(n)}}_{[0,T]}} E_n
H^{S^{(n)}}(\sigma,\Psi_n(\tau))=V_n=\inf_{\sigma\in\mathcal
{T}^{S^{(n)}}_{[0,T]}}\sup_{\tau\in
\Delta_n} E_n H^{S^{(n)}}(\Phi_n(\sigma),\tau).
\]
Thus
%
\begin{eqnarray}\label{5.13}
\Gamma\bigl(S^{(n)}\bigr)-V_n&\leq&\inf_{\sigma\in\Delta_n}\sup_{\tau\in
\mathcal{T}^{S^{(n)}}_{[0,T]}} E_n H^{S^{(n)}}(\sigma,\tau)\nonumber\\
&&{}- \inf_{\sigma\in\Delta_n}\sup_{\tau\in
\mathcal{T}^{S^{(n)}}_{[0,T]}} E_n
H^{S^{(n)}}(\sigma,\Psi_n(\tau))\nonumber\\
&\leq&\sup_{\tau\in
\mathcal{T}^{S^{(n)}}_{[0,T]}} E_n\bigl|Y^{S^{(n)}}_{\tau}-
Y^{S^{(n)}}_{\Psi_n(\tau)}\bigr|\quad\mbox{and} \\
V_n-\Gamma\bigl(S^{(n)}\bigr)&\leq&\inf_{\sigma\in
\mathcal{T}^{S^{(n)}}_{[0,T]}}\sup_{\tau\in\Delta_n} E_n
H^{S^{(n)}}(\Phi_n(\sigma),\tau)\nonumber\\
&&{}- \inf_{\sigma\in\mathcal{T}^{S^{(n)}}_{[0,T]}}\sup_{\tau\in
\Delta_n} E_n H^{S^{(n)}}(\sigma,\tau)\nonumber\\
&\leq&\sup_{\sigma\in
\mathcal{T}^{S^{(n)}}_{[0,T]}} E_n\bigl|X^{S^{(n)}}_{\Phi_n(\sigma)}-
X^{S^{(n)}}_{\sigma}\bigr|.\nonumber
\end{eqnarray}
For any $0\leq t_1,t_2\leq{T}$, and $u_1,u_2\geq{0}$, we have the
following inequalities: $|e^{-rt_1}u_1-e^{-rt_2}u_2|\leq
|u_2-u_1|+r|t_1-t_2|\max(u_1,u_2)$ and
$|e^{rt_1}u_1-e^{rt_2}u_2|\leq
e^{rT}(|u_2-u_1|+r|t_1-t_2|\max(u_1,u_2))$. This together with
(\ref{5.13}) gives
%
\begin{eqnarray}\label{5.14}
&&\bigl|V_n-\Gamma\bigl(S^{(n)}\bigr)\bigr|\nonumber\\
&&\qquad\leq\sup_{\tau\in
\mathcal{T}^{S^{(n)}}_{[0,T]}} E_n\bigl(\bigl|Y^{S^{(n)}}_{\tau}-
Y^{S^{(n)}}_{\Psi_n(\tau)}\bigr|+\bigl|X^{S^{(n)}}_{\Phi_n(\tau)}
- X^{S^{(n)}}_{\tau}\bigr|\bigr) \nonumber\\[-8pt]\\[-8pt]
&&\qquad\leq2\frac{rT}{n}\Bigl(M+E_n\sup_{0\leq t \leq
T} e^{rt}S^{(n)}_{t}\Bigr)+2e^{rT}\frac{rT}{n}E_n\sup_{0\leq t \leq T}
S^{(n)}_{t}\nonumber\\
&&\qquad\quad{}+(\delta+2e^{rT})\sup_{\tau\in
\mathcal{T}^{S^{(n)}}_{[0,T]}}E_n\bigl|S^{(n)}_{T\wedge(\Phi_n(\tau)+
{T}/{n})}-S^{(n)}_{\Phi_n(\tau)}\bigr|.\nonumber
\end{eqnarray}
The sequence $S^{(n)}$, $n\geq1$, satisfies Assumption
\ref{asm2.3}, and so from (\ref{5.14})
%
\begin{eqnarray}\label{5.15}
&&\lim_{n\rightarrow\infty}\bigl|V_n-\Gamma\bigl(S^{(n)}\bigr)\bigr|\nonumber\\
&&\qquad\leq
(\delta+2e^{rT})\lim_{n\rightarrow\infty}\sup_{\tau\in
\mathcal{T}^{S^{(n)}}_{[0,T]}}E_n\bigl|S^{(n)}_{T\wedge(\Phi_n(\tau)+
{T}/{n})}-S^{(n)}_{\Phi_n(\tau)}\bigr|\\
&&\qquad=0.\nonumber
\end{eqnarray}
\upqed\end{pf}
\begin{rem}\label{rem5.1}
Similar results can be obtained for game options in the Merton model
with integral payoffs (Asian options). It can be checked that
integral payoffs satisfy Assumption \ref{asm2.1} and by estimates in
the spirit of this section we can get convergence results for this
case also. Of course, put and call options can be treated even in a
more simple way since their payoffs depend only on the present stock
price.
\end{rem}

\section*{Acknowledgments}
I am deeply grateful to my adviser and
teacher, Yuri Kifer, for guiding me and helping me to present this
work. I also want to thank the referees for valuable
suggestions.\vadjust{\goodbreak}

%

%
\printaddresses

\end{document}